\documentclass[preprint,authoryear,12pt]{elsarticle}

\usepackage[fleqn]{amsmath}
\usepackage[]{amsfonts}
\usepackage[overload]{empheq}
\usepackage{color}
\newcommand{\R}{\mathbb{R}}
\journal{Transportation Research Part B}

\renewcommand{\P}{\mathcal{V}_{BP}}

\renewcommand{\L}{\mathcal{L}}

\newcommand{\BP}{\Gamma^P}
\newcommand{\BPP}{\mathcal{Y}}
\newcommand{\B}{\Gamma}
\newcommand{\D}[1]{\Delta #1}

\renewcommand{\and}{\hskip 2mm \mbox{and} \hskip 2mm}

\begin{document}

\begin{frontmatter}

\title{The impact of source terms in the variational representation \\ of traffic flow}

\author[GATECH]{Jorge A. Laval\corref{cor}} 
\ead{jorge.laval@ce.gatech.edu}
\cortext[cor]{Corresponding author. Tel. : +1 (404) 894-2360; Fax : +1 (404) 894-2278}
\author[GUIL]{Guillaume Costeseque} 
\author[GATECH]{Bargavarama Chilukuri} 
\address[GATECH]{School of Civil and Environmental Engineering, Georgia Institute of Technology}
\address[GUIL]{Inria Sophia Antipolis - M\'{e}diterran\'{e}e, France.}

\begin{abstract}

This paper revisits the variational theory of traffic flow, now under the presence of continuum lateral inflows and outflows to the freeway
{say Eulerian source terms}. 
It is found that a VT solution exists only  in Eulerian coordinates when source terms are exogenous
{meaning that they only depend on time and space}, 
but not when they are a function of traffic conditions, as per a merge model. In discrete time, however, these dependencies become exogenous, which allowed us to propose improved numerical solution methods. In Lagrangian and vehicle number-space coordinates, VT solutions 
{may} 
not exist even if source terms are exogenous.

\end{abstract}


\end{frontmatter}

%\linenumbers           % Running line numbers.

\section{Introduction}

The variational theory (VT) 
{applied in}
traffic flow
{theory} 
\citep{Dag05a, Dag05b} was an important milestone. Previously, the only analytical solution to the kinematic wave model of \cite{Lig55,Ric56} was the method of characteristics, which does not give a global solution
{in time}
as one needs to keep track of characteristics crossings (shocks
{and rarefaction waves} )
and impose entropy conditions to ensure uniqueness. This means that analytical solutions cannot be formulated except for very simple problems.

In contrast, VT makes use of the link between conservation laws and the Hamilton-Jacobi partial differential equation (HJ PDE), allowing the kinematic wave model to be solved using the Hopf-Lax formula \citep{Lax1957, Olejnik1957, Hopf1970}, better known in transportation as Newell's minimum principle \cite{New93}
{whenever the fundamental diagram or Hamiltonian is piecewise linear}. 
The big advantage is that this
{representation}
formula gives an analytical global solution
{in time}
that does not require explicit consideration of shocks and/or entropy conditions. Moreover, current approximation methods for the Macroscopic Fundamental Diagram (MFD) of urban networks rely on this approach \citep{Daganzo08a, Geroliminis13a, leclercq2013estimating, Lav15}. In addition, when the flow-density fundamental diagram is triangular
{(or piecewise linear more generally)},
numerical solutions become exact in the HJ  framework, and this is not the case  in the conservation law approach \citep{Lev93}. It has also been shown that the traffic flow problem cast  in Lagrangian or vehicle number-space coordinates also accepts a VT solution, which can be used to obtain very efficient numerical solution methods \cite{Lec07,Lav13hj}.

The traffic flow problem with source term is also of great interest; e.g., it can be used to approximate (i)  long freeways with closely spaced entrances and exits, (ii) the effect of lane-changing activity on a single lane, or (iii) the effects of turning movements, trip generation and trip ends in the MFD. But the underlying assumption in the previous paragraphs is that vehicles are conserved. This begs the
{twofold}
question, are
{representation formulas for the}
VT solutions still valid, or even applicable, when there is a
{Eulerian}
source term? If not, can efficient numerical solution methods still be implemented? Recent developments in this area have not answered these questions as they are primarily concerned with discrete source terms \citep[e.g.,][]{Dag14,Cos14,costeseque2014variational}.

To answer these questions this paper is organized as follows.  In section~\ref{sec:pbm} we formulate the general problem and show that in general VT solutions are not applicable; but section~\ref{sec:exogenous} shows that they are when the source term is exogenous. Based on these results, section~\ref{sec:numerical} presents numerical methods for the endogenous inflow problem that outperform existing methods. Section~\ref{sec:othercoord} briefly shows that in space-Lagrangian and time-Lagrangian  coordinates VT solutions do not exist even if source terms are exogenous. A discussion of results and outlook is presented in section~\ref{sec:disc}.

\section{Problem Formulation}
\label{sec:pbm} 

Consider a long homogeneous freeway corridor with a large number of entrances and exits such that the net lateral freeway inflow rate, $\phi$, or inflow for short, can be treated as a continuum variable in time $t { \geq 0}$ and location $x {\in \R}$, and has units of veh/time-distance. The inflow is an endogenous variable consequence of the demand for travel, and could be captured by  a function of the traffic states both in the freeway and the ramps. For simplicity, in this paper the inflow is assumed to be a function of the density, $k(t,x)$, of the freeway only, i.e. $\phi=\phi(k(t,x))$, but also the exogenous case $\phi=\phi(t,x)$ will be of interest. 

In any case, the traffic flow problem analyzed in this paper is the following conservation law with source term:
\begin{subequations}\label{k-model}
\begin{align}
	k_t+H(k)_x &=\phi,\label{a}\\
	k&=g \hskip .5cm \mbox{on} \hskip .3cm \Gamma\label{b}
\end{align}
\end{subequations}
where $H$ is the fundamental diagram, $g$ is the data defined on a boundary $\Gamma$, and variables in subscript represent partial derivatives. Now we define the function $	N(t,x)$ such that:
\begin{equation}\label{N}
N_x=- \ k 
{\quad \mbox{say} \quad N(t,x) := \int_{x}^{+\infty} k(t,y) dy},
\end{equation}
and integrate (\ref{k-model}) with respect to $x$ to obtain its HJ form \citep{Eva98}:
\begin{subequations}\label{N-model}
\begin{align}
	N_t-H(-N_x)&=\Phi,\label{N-modela}\\
	N&=G \hskip .5cm \mbox{on} \hskip .3cm \Gamma,\label{N-modelb}
\end{align}
\end{subequations}
where we have defined:
\begin{subequations}\label{GPhi}
\begin{align}
	G(t,x)&=\oint _\Gamma g(t,x) d\Gamma, \hskip .5cm (t,x)\in\Gamma, \hskip .5cm \mbox{and}\label{G}\\
	\Phi(t,x)&=-\int_{0}^x\phi (t,y) dy \label{Phi}
\end{align}
\end{subequations}
where $\Phi$ is a potential function; the negative sign in  (\ref{Phi}) follows from the traditional counting convention in traffic flow, where further downstream vehicles have lower vehicle numbers. 

It is important to note that compared to the traditional case with zero inflow, under a continuum source term the interpretation of $N$ changes: for fixed $x$ it
{still defines a cumulative count curve due to~\eqref{N}
but}
its isometrics do no longer give vehicle trajectories. 
To see this we note that according to (\ref{N-modela}) the flow $q=H(k)$ is now:
\begin{equation}\label{flow}
	q=N_t-\Phi,
\end{equation}
The time integration of (\ref{flow}) reveals that cumulative count curves are now given, up to an arbitrary constant, by $\tilde{N}(t,x)+\int_0^t\int_{0}^x\phi(s,y)dsdy$, where
{$\tilde{N}(t,x) := \int_{0}^{t} q(s,x) ds$ denotes the usual $N$-curve
without source terms}
and the integral represents the net number of vehicles entering
{($\phi > 0$) or exiting ($\phi < 0$)}
the road segment by time $t$ and upstream of $x$.

The method to obtain the solution of (\ref{N-model}) depends on the dependencies of the potential function. If the inflow function is allowed to depend on the traffic state in the freeway, i.e. $\phi=\phi(k)$ then the potential function $\Phi$ depends non-locally on $N$, since:
\begin{equation}\label{}
	\Phi (t,x) =\tilde{\Phi} (N,x)=-\int_{0}^x\phi(-N_x (t,y)) dy ,
\end{equation}
and therefore
{it is not obvious that (\ref{N-model}) accepts a representation formula
as a VT solution}. 
To see this, we note that even in the simplest linear case:
\begin{subequations}\label{ephi}
\begin{align}
	\phi(k)&=a-bk, \hskip 1cm a,b \geq 0, \hskip .5cm \mbox{we get:}\label{phia}\\
	\tilde{\Phi} (N,x)&=-ax+\left(\int_0^x k(t,y) \ dy \right) b =-ax-( N(t,x) -c ) b,\label{Phiendo2}
\end{align}
\end{subequations}
where $c {= \int_{0}^{+\infty} k(t,y) dy = N(t,0)}$ 
is an arbitrary constant of integration. 
This means that (\ref{N-model}) becomes the more general HJ equation $N_t-\tilde{H}(x,N,-N_x)=0$, where $\tilde{H}$ is the Hamiltonian. The Hamiltonian's $N$-dependency is what complicates matters. \cite{Barron1996,barron2015representation} show that a Hopf-Lax type solution exists in this case only when 
{(among other assumptions)}
$a=0$ 
{i.e. $\tilde{H} (x, u, p) = \hat{H} (u, p)$ 
where $\hat{H}$ does not depend on the space variable}
and  $\tilde{H}$ is homogeneous of degree one 
{with respect to} 
$N_x$, which is of no use in traffic flow
{since in practice it means that the Hamiltonian has to be monotone linear}. 
As explained in section 5, we argue that this $N$-dependency prevents formulating an equivalent VT problem (such as  (\ref{Nsol}) below) that can be solved with variational methods. 

As shown in the next section, it turns out that when  inflows are exogenous  a global VT solution can still be identified, albeit not in Hopf-Lax form because optimal paths are no longer straight lines. 

\section{Exogenous inflow}
\label{sec:exogenous}

\subsection{Setting of the VT problem}

\begin{figure}[tb]
\begin{center}
\includegraphics[width=\textwidth]{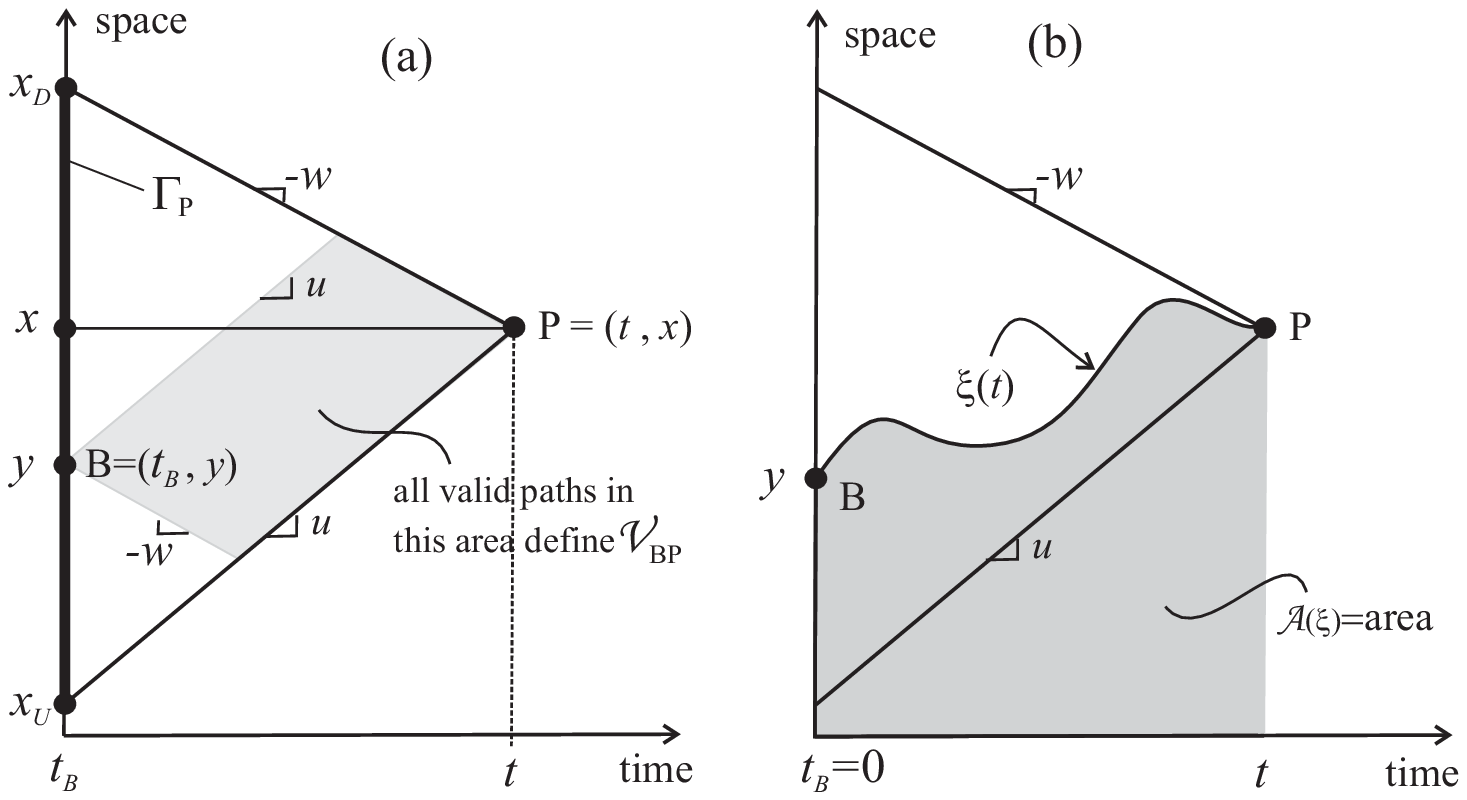}
\end{center}
\caption[]{Illustration of key definitions in VT: (a) the initial value problem; (b) the area of the integration to obtain the $J$-integral.}
\label{f1}
\end{figure}

The results in this section are based on VT \citep{Dag05a}, where the solution of the HJ equation $N_t-\tilde{H}(t,x,-N_x)=0$ is given by the solution of the following variational problem: 
\begin{subequations}\label{Nsol}
\begin{align}
	N(P) &= \min_{B\in\BP,\ \xi\in\P}f(B,\xi), \hskip .5cm \mbox{with}\label{Nsola}\\
	f(B,\xi)&=G(B) + \int_{t_B}^t R(s,\xi(s),\xi'(s)) \ ds  \label{Nsolb}\end{align}
\end{subequations}
where $P$ is a generic point  with coordinates $(t,x)$, $B\equiv(t_B,y)$ is a point in the boundary $\BP$, $\xi$ is a member of the set of all valid paths between  $B$ and $P$ denoted $\P$, and $\xi(t_B)=y$; see Fig.~\ref{f1}a. The function $R(\cdot)$ gives the maximum passing rates along the observer and corresponds to the {(concave)} Legendre transform of $\tilde{H}$, i.e.,
$$ R(t,x,v)=\sup_{k} \left\{ \tilde{H}(t,x,k)-vk \right\}. $$

It is worth mentioning that in the simplest homogeneous case where $\tilde{H}=H(k)$ the VT solution (\ref{Nsol}) becomes the Hopf-Lax formula:
\begin{equation}\label{HL}
N(P) = \min_{B\in\BP}\left\{G(B)+(t-t_B)R\left(\frac{x-y}{t-t_B}\right)\right\},
\end{equation}
where the minimization over $\xi(t)$ is no longer necessary since characteristics become straight lines in this case. Unfortunately, even the presence of exogenous lateral inflows that vary in time or space make characteristics not to be straight lines and therefore a Hopf-Lax type  solution cannot be devised. A VT solution, however, still exists as shown next. \\

We now formulate VT solution to account for source terms explicitly.
{In the remaining of the paper}, 
we assume a triangular flow-density diagram. It may be defined by its free-flow speed $u$, wave speed $-w$ 
{(with $w > 0$)}
and jam density $\kappa$ 
{such that
\begin{equation}\label{eq:triangFD}
	H(k) = \min \left\{ u k \ , \ w (\kappa - k) \right\} 
	\quad \mbox{for any} \quad k \in [0, \kappa].
\end{equation} }
It follows that the capacity is $Q=\kappa w u /(w + u)$ and the critical density $K=Q/u$. \\

When the potential function is exogenous, i.e. $\Phi=\Phi(t,x)$, the Hamiltonian can be written as the sum of the fundamental diagram and the potential function, i.e.: $\tilde{H}(t,x,k)=H(k)+\Phi(t,x)$. In the case of a triangular fundamental diagram we have $R(t,x,v)=Q-Kv+\Phi(t,x)$
{with $v \in [-w, u]$},
and the function $f(B,\xi)$ to be minimized reads:
\begin{equation}\label{fxB00}
	f(B,\xi)=G(B)  + (t-t_B)Q-(x-y)K +\underbrace{\int_{t_B}^t \Phi(s,\xi(s)) \ ds}_J
\end{equation}
The $J$-integral in (\ref{fxB00}) is what separates this problem from the problems studied so far in traffic flow using VT principles. This integral in terms of $\phi$ is:
\begin{equation}\label{J}
J=- \int_{t_B}^t\int_{y}^{\xi(s)} \phi(s,x) \ dx ds, 
\end{equation}
and represents the net number of vehicles \emph{leaving} the area below the curve $x=\xi(t)$, namely area $A(\xi)$; see  Fig.~\ref{f1}b. Therefore, minimizing $J$ given $y=\xi(t_B)$, can be interpreted as finding $\xi(t)$ that maximizes the net number of vehicles \emph{entering} $A(\xi)$.

\subsection{Initial value problems}

In the initial value problem (IVP) the boundary $\B$ is the line $\left\{ t_B=0 \right\} {\times \R}$,  so that the problem here is (\ref{N-modela}) supplemented with:
\begin{equation}\label{CIDCI}
	N(0,x)=G(x),
\end{equation}
{where we assume that $G \in C^2 (\R)$.}
The candidate set for $B, \BP$ is reduced to $B$'s $x$-coordinate, $y$, which is delimited by two points $U=(0,x_U)$ and $D=(0,x_D)$, where:
\begin{subequations}\label{xuxd}
\begin{align}
	x_U&=x-ut,\and	x_D=x+wt,\label{b}\\
	x_U&<y<x_D\label{ycond}
\end{align}
\end{subequations}
see Fig.~\ref{f1}a. The following subsections examine simplified versions of this problem that reveal  considerable insight into the general solution.

\subsubsection{Constant inflow}

Consider the IVP with constant inflow problem:
\begin{equation}\label{}
	\phi(t,x)=a,
\end{equation}
{for some $a \in \R \setminus \{ 0 \}$.}
It follows from  (\ref{J}) that $J=-aA(\xi)$ and therefore, for a fixed $y=\xi(0)$ the minimum of  (\ref{fxB00}) is obtained by a path that: (i) maximizes $A(\xi)$ when $a>0$; or (ii) minimizes $A(\xi)$ when $a<0$.
These two optimum paths are the extreme paths in  $\P$ that define its boundary; see  ``upper" and  ``lower" paths in Fig.~\ref{f2}a. This solution is a ``bang-bang" solution, typical for this type optimal control problems. The reader can verify that the  areas under these paths are given by $A(y,x_D)$ and $A(y,x_U)$, respectively, where we have defined:
\begin{equation}\label{}
	A(y,x_-)=\frac{1}{2} \left((x+x_-)t- \text{sign}(a) \frac{\left(x_- -y \right){}^2}{u+w}\right)
\end{equation}
where $\text{sign}(a)$ is the sign of $a$ and $x_-$ is a placeholder for $x_D$ if $a>0$ or $x_U$ if $a<0$. Thus the optimization problem has been reduced to a single variable, $y$, whose first- and second-order conditions for a minimum {$y^*$} read:
\begin{subequations}\label{}
\begin{align}
	f'(y^*)&= \psi \left(y^*-x_-\right)+G'(y^*)+K=0,\label{a}\\
	f''(y^*)&= \psi+G''(y^*)>0, \hskip .5cm \mbox{where:}\label{second_order}\\
	\psi& :=\frac{ |a| }{u+w},\label{psi}
\end{align}
\end{subequations}
{where $| \cdot |$ denotes the absolute value. 
Notice that $\psi > 0$ as long as $a \neq 0$.}

 Of course, the optimal point $y^*$ needs to satisfy  (\ref{ycond}):  $x_U\le y^* \le x_D$.

\subsubsection{Constant initial density}
Here, in addition to $\phi(t,x)=a,$ we assume
\begin{equation}\label{}
	g(x)=k_0,\hskip 1cm -\infty<x<\infty 
\end{equation}
which implies that the function to minimize is a parabola:
\begin{equation}\label{}
	f(y)= - c_0 - c_1 y + \frac{\psi}{2} y^2,
\end{equation}
with constants  $c_0=K x_U - \dfrac{\psi}{2} \left[ x_-^2 { - \text{sign}(a) } \left(x_- +x \right) (u+w)t \right]$ and $c_1= \psi x_- -(K-k_0)$, and extremum:
\begin{equation}\label{xbt}
	y^*=x_- - \dfrac{ K-k_0 }{\psi} .
\end{equation}

Notice that (\ref{second_order}) is always satisfied in this case and therefore $y^*$ is always a minimum and should be included so long as $x_U\le y^*\le x_D$. Combining this with  (\ref{xbt}) gives that the final solution can be expressed as:
\begin{subequations}\label{NconstSol}
\begin{empheq}[left={N(t,x) =\empheqlbrace}]{align}
     &f(y^*), &t>(K-k_0)/a>0\label{NconstSola}\\
     &\min\{f(x_U),f(x_D)\}, &\mbox{otherwise}
\end{empheq}
\end{subequations}
Notice that $(K-k_0)/a$ represents the time it takes for the system to reach critical density, namely ``time-to-capacity". This means that $y^*$ will be the optimal candidate only once a regime transition occurs. This will happen only if the time to capacity is positive, i.e. when sign($a$) = sign($K-k_0)$
{
or more explicitly, when $k_0$ is under-critical and $\phi$ is an inflow ($a$ positive)
or when $k_0$ is over-critical and $\phi$ is an outflow ($a$ negative). 
}

Somewhat unexpectedly, we note that $-f_x(x_U)=-f_x(x_D)=-f_x(y^*)=k_0+at$, which means that the density is always given by {the traveling wave} 
$ k(t,x)=k_0+at $,
which is also the solution of  (\ref{k-model}) in this case using the method of characteristics. One should also impose feasibility conditions for the density, i.e. $0\le k\le \kappa$, which gives:
\begin{subequations}\label{kconstSol}
\begin{empheq}[left={k(t,x) =\empheqlbrace}]{align}
     &0 , &t>-k_0/a, a<0\label{kconstSola}\\
     &\kappa , &t>(\kappa-k_0)/a, a>0\label{kconstSola}\\
     & k_0+at, &\mbox{otherwise}
\end{empheq}
\end{subequations}
The flow can be obtained using (\ref{flow}), but it is equivalent and simpler to use $q=H(k) $. 

\subsubsection{Extended Riemann problems}

\begin{figure}[tb]
\begin{center}
\includegraphics[width=\textwidth]{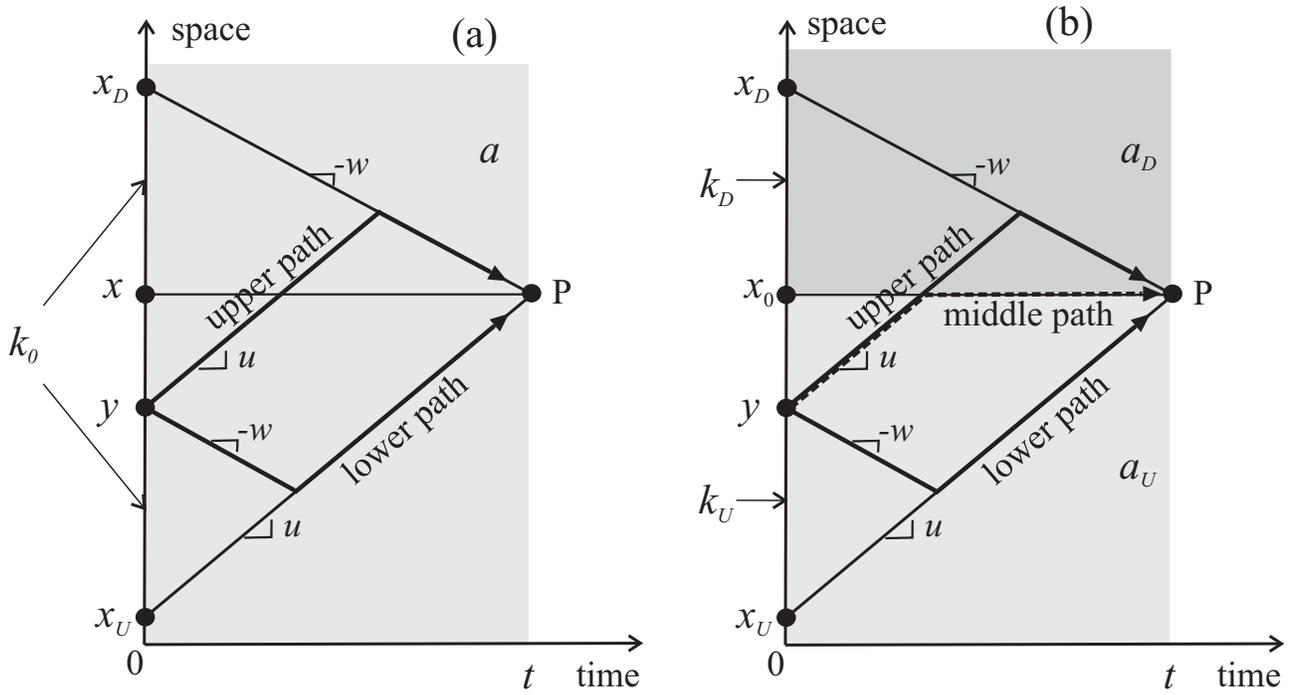}
\end{center}
\caption[]{Possible paths to minimize the $J$-integral: (a) constant initial density; (b) extended Riemann problems.}
\label{f2}
\end{figure}

Riemann problems are the building blocks of Godunov-type numerical solution methods. As illustrated in Fig.~\ref{f2}b, in these problems one is interested in the value of $N$ at $x=x_0\ge ut$, i.e. at point  $P=(t,x_0)$, with initial data typically given by the density at $t=0$. 
{
This is a very special case for with the `target' point $P$ is located
on the discontinuity of the initial data at $x=x_0$.
However the methodology below can be easily extended to the cases
where $P=(t,x)$ with $x < x_0$ or $x > x_0$.
}
Here, we extend the initial data to include the inflow:
\begin{subequations}
\begin{empheq}[left={(g(x),\phi(x)) =\empheqlbrace}]{align}
     &(k_U,a_U), \hskip 1cm x\le x_0\\
     &(k_D,a_D), \hskip 1cm x>x_0,
\end{empheq}
\end{subequations}
which in conjunction with  (\ref{N-model}) define an extended Riemann problem, or ERP for short.
{We assume that $\left( a_U, a_D \right) \neq (0,0)$.}
Notice that now:
\begin{equation}\label{xuxd}
	x_U=x_0-ut,\and	x_D=x_0+wt.
\end{equation}

For simplicity and without loss of generality we set $G(x_0)=0$, which implies:
\begin{subequations}
\begin{empheq}[left={G(x) =\empheqlbrace}]{align}
     &(x_0-x)k_U, \hskip 1cm x\le x_0\\
     &(x_0-x)k_D, \hskip 1cm x>x_0
\end{empheq}
\end{subequations}
It will be convenient to define:
\begin{equation}\label{}
	\eta = a_U/a_D,\hskip 1cm	\psi=\frac{a_D}{u+w},\hskip 1cm \theta=u/w.
\end{equation}
The $J$-integral in this case is a weighted average of the portion of $A(\xi)$ upstream and downstream of $x=x_0$, weighted by $a_U$ and $a_D$, respectively. For instance, if $\xi(t)>0$ for all $t$, we have 
$
J= -\int_{0}^t\int_{0}^{\xi(t)} \phi(x) \ dx dt=-\int_{0}^t(a_Ux_0+a_D(\xi(t)-x_0)) dt=- (a_Ux_0t+a_D\int_{0}^t\xi(t)-x_0 dt), 
$
where $x_0t$ is the portion of $A(\xi)$ upstream of $x=x_0$ and  $\int_{0}^t\xi(t)-x_0 dt$ is the portion of $A (\xi)$ downstream of $x=x_0$ in this particular case.

It follows that the minimization of the $J(\xi)$ can be achieved analogously to the previous section by considering the upper and lower paths from each candidate $y=\xi(0)$. In addition, however, one has to include 2 middle ``paths" that would reach and stay at $x=x_0$ until reaching $P$; see  Fig.~\ref{f2}b. To formalize, let $j_1, j_2$ and $j_3$ be the value of $J(\xi)$ when $y\ge x_0$ along the upper, lower and middle paths, respectively; similarly for $j_4, j_5$ and $j_6$ when $y\le x_0$. Calculating the areas upstream and downstream of $x=x_0$ defined by each   path, it can be shown that the  $j_i$'s can be obtained from:
\begin{subequations}\label{ji}
\begin{align}
	2(j_1(y)-J_0)/\psi &=  -\theta  t^2 w^2 (2 \eta  (\theta +1)+1)+2 t w y_0+y_0^2   ,\\
	2(j_2(y)-J_0)/\psi &=- \eta  \theta  (2 \theta +1) t^2 w^2+2 \eta  \theta  t w y_0+y_0^2 ((\eta -1) \theta -1)   ,\\
	2(j_3(y)-J_0)/\psi &= -(\theta +1) \left(2 \eta  \theta  t^2 w^2+y_0^2\right)   ,\\
	2(j_4(y)-J_0)/\psi &= -\theta  t^2 w^2 (2 \eta  (\theta +1)+1)+2 t w y_0+y_0^2 (\eta  \theta +\eta -1)/\theta   ,\\
	2(j_5(y)-J_0)/\psi &= -\eta  \left(\theta  (2 \theta +1) t^2 w^2-2 \theta  t w y_0+y_0^2\right)   ,\\
	2(j_6(y)-J_0)/\psi &= \eta  (\theta +1) \left(y_0^2-2 \theta ^2 t^2 w^2\right)/\theta  ,
\end{align}
\end{subequations}
where $y_0=x_0-y$ and $J_0$ is a constant of our problem that represents the net number of vehicles \emph{leaving} the area upstream of $x_U$, i.e. $-a_Ux_Ut$ in this case. Notice that $j_1(x_0)=j_4(x_0), j_2(x_0)=j_5(x_0)$ and $j_3(x_0)=j_6(x_0)$, as expected. The function to minimize in this case can be written as:
\begin{equation}\label{}
	f(y)=G(y) + tQ-(x_0-y)K + J(y), \hskip .5cm \mbox{where:}
\end{equation}
\begin{subequations}\label{JJ}
\begin{empheq}[left={J(y) =\empheqlbrace}]{align}
     &\min\{j_1(y),j_2(y),j_3(y)\}, &y>x_0\label{2nda}\\
     &\min\{j_4(y),j_5(y),j_6(y)\}, &y\le x_0\label{2nda}
\end{empheq}
\end{subequations}

The solution of (\ref{N-model}) under these conditions can be reduced to the evaluation of $f(y)$ at a small number of candidates. In addition to candidates $y=x_U$ and $y=x_D$, we have to consider the discontinuity at  $y=x_0$ and the possible minima produced by each of the components in  (\ref{ji}). Since $f(y)$ is piecewise quadratic, each one of these components has at most one minimum, namely  $y=y_i, i=1,\ldots 6$, which can be obtained by solving the first-order conditions $f(y_i)=0$ associated with each of the $j_i$'s; i.e.:
\begin{subequations}\label{}
\begin{align}
	y_1&=x_D-(K-k_D)/\psi ,\hskip 1cm 
	y_2=x_0+\frac{\eta  \theta  x_D-(K-k_D)/\psi }{ \theta(\eta-1)-1 },\label{b}\\
	y_3&=x_0+\frac{K-k_D}{(\theta+1) \psi } ,\hskip 1.6cm 
	y_4=x_0+\frac{x_U+\theta(K-k_U)/\psi  }{ 1-\eta( 1+ \theta)   },\label{b}\\
	y_5&=x_U+(K-k_U)/(\eta\psi),\hskip .5cm 
	y_6=x_0+\frac{\theta  \left(k_U-K\right)}{\eta  (\theta +1) \psi },\label{b}
\end{align}
\end{subequations}
For the $y_i$'s to be valid candidates they must meet the following conditions:
\begin{subequations}\label{yicond}
\begin{align}
	x_0<y_i<x_D, \hskip 2cm &i=1,2,3\label{a}\\
	x_U<y_i<x_0, \hskip 2cm &i=4,5,6\label{a}
\end{align}
\end{subequations}
which ensure that $y_i$ is on the boundary $\BP$. With all, the sought solution can be expressed as:
\begin{equation}\label{Nsolriem}
	N(t,x_0)= \min_{y\in\BPP}f(y), \hskip .5cm \mbox{with: }
	\BPP=\{x_U,x_0,x_D,y_1^*,\ldots y_6^*\}
\end{equation}
where $y_i^*$ is $y_i$ if (\ref{yicond}) is met and null otherwise. Notice that this solution method does not require imposing $f''(y^*)>0$ because maxima will be automatically discarded in the minimum operation.

The average flow $\bar{q}(t,x_0)$ during $(0,t)$ can be obtained as follows:
\begin{subequations}\label{}
\begin{align}
	\bar{q}(t,x_0)&=\frac{1}{t}\int_0^t N_t(s,x_0)-\Phi(s,x_0) \ ds
	=\frac{1}{t}(N(t,x_0)-\int_0^t \Phi(s,x_0) \ ds)\label{b}\\
	&=\frac{1}{t}(N(t,x_0)-(J_0+x_Ua_U)t) =N(t,x_0)/t-J_0-x_Ua_U .\label{avgflow}
\end{align}
\end{subequations}

As shown in Section~\ref{sec:numerical}, this formula allows to formulate improved numerical solution methods for the endogenous problem.

{
\subsection{Initial and boundary value problem}
}

{
We restrict the space to a segment of road $[\xi,\chi]$ with $\xi < \chi$.
In the initial and boundary value problem (IBVP) the boundary $\B$ is the set defined as 
$$\left( \left\{ t_B=0 \right\} \times [\xi, \chi] \right) \cup 
\left( (0,+\infty) \times \left\{ x=\xi \right\} \right) \cup
\left( (0,+\infty) \times \left\{ x=\chi \right\} \right), $$ 
so that the problem here is (\ref{N-modela}) supplemented with:
\begin{equation}\label{CIBVP}
\begin{cases}
	N(0,x)=G_{ini}(x), \quad &\mbox{on} \quad [\xi,\chi], \\
	N(t,\xi)=G_{up}(t), \quad &\mbox{on} \quad (0,+\infty), \\
	N(t,\chi)=G_{down}(t), \quad &\mbox{on} \quad (0,+\infty).
\end{cases}
\end{equation}
}

{
For obvious compatibility reasons, we request these conditions to satisfy
$$ G_{ini} (\xi ) = G_{up}(0) \quad \mbox{and} \quad G_{ini}(\chi) = G_{down}(0).$$
We also consider a $a \neq 0$ and $k_0 \in [0, \kappa]$ 
such that the inflow rate is given by
$$ \varphi (t,x) = a, \quad \mbox{for any} \quad (t,x) \in [0, +\infty) \times [\xi, \chi]$$
and the initial density
$$ g_{ini} (x) = k_0 , \quad \mbox{for any} \quad x \in [\xi, \chi].$$}

\begin{figure}[tb]
\begin{center}
\includegraphics[width=.7\textwidth]{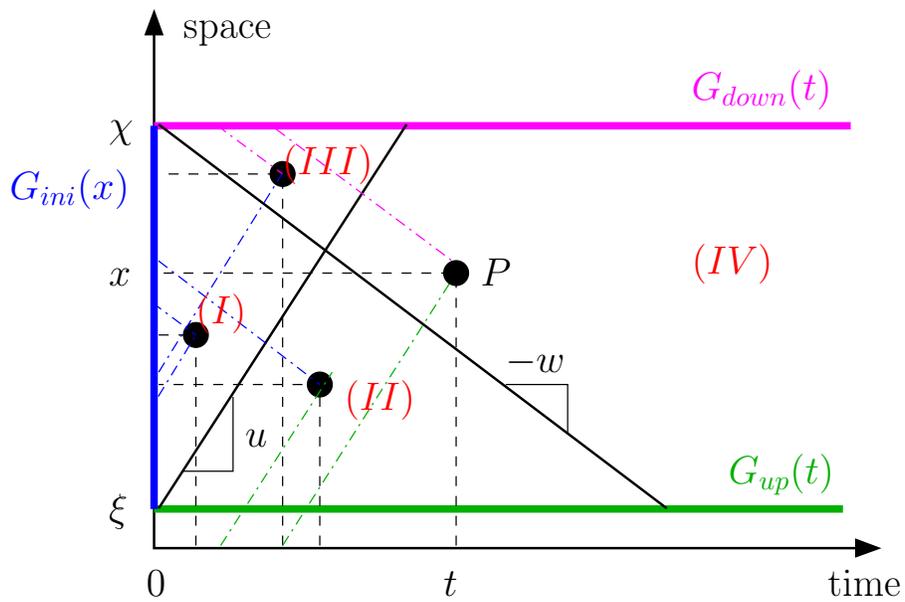}
\end{center}
\caption[]{The four different regions to consider for an Initial and Boundary Value Problem with initial data $G_{ini}$, upstream $G_{up}$ and downstream $G_{down}$ boundary conditions.}
\label{f6}
\end{figure}

{
According to the position of point $P=(t,x)$ with $t > 0$ and $x \in (\xi, \chi)$, 
we have the following cases to distinguish
(in the spirit of \cite{jin2015continuous})
$$ \begin{aligned}
&\left( t_U, x_U \right) = \begin{cases}
\left( 0, x - \dfrac{t}{u} \right), \quad &\mbox{if} \quad x \geq \xi + \dfrac{t}{u},  \\
\left( t - \dfrac{x- \xi}{u}, \xi \right), \quad &\mbox{else},
\end{cases}
\\
\mbox{and}& \\
&\left( t_D, x_D \right) = \begin{cases}
\left( 0, x + \dfrac{t}{w} \right), \quad &\mbox{if} \quad x \leq \chi - \dfrac{t}{w},  \\
\left( t + \dfrac{x + \chi}{w}, \chi \right), \quad &\mbox{else},
\end{cases}
\end{aligned}$$
that define 4 regions (see Figure~\ref{f6})
{\footnotesize
$$\begin{array}{rccc}
\text{Region} & \text{Set of points} & \text{Upstream point} & \text{Downstream point} \\
I: & \left\{ (t,x) \ \left| \ \chi - \dfrac{t}{w} \geq x \geq \xi + \dfrac{t}{u} \right. \right\} & \left( t_U, x_U \right) = \left( 0, x - \dfrac{t}{u} \right) & \left( t_D, x_D \right) = \left( 0, x + \dfrac{t}{w} \right) \\
\\
II: & \left\{ (t,x) \ \left| \ x \leq \min \left\{ \xi + \dfrac{t}{u}, \chi - \dfrac{t}{w} \right\} \right. \right\} & \left( t_U, x_U \right) = \left( t - \dfrac{x- \xi}{u}, \xi \right) & \left( t_D, x_D \right) = \left( 0, x + \dfrac{t}{w} \right) \\
\\
III: & \left\{ (t,x) \ \left| \ x \geq \max \left\{ \xi + \dfrac{t}{u}, \chi - \dfrac{t}{w} \right\} \right. \right\} & \left( t_U, x_U \right) = \left( 0, x - \dfrac{t}{u} \right) & \left( t_D, x_D \right) = \left( t + \dfrac{x + \chi}{w}, \chi \right) \\
\\
IV: & \left\{ (t,x) \ \left| \ \xi + \dfrac{t}{u} \geq x \geq \chi - \dfrac{t}{w} \right. \right\} & \left( t_U, x_U \right) = \left( t - \dfrac{x- \xi}{u}, \xi \right) & \left( t_D, x_D \right) = \left( t + \dfrac{x + \chi}{w}, \chi \right)
\end{array} $$
}
}

{
To be continued 
but the idea is that in Region I, everything is exactly the same than in previous Subsection 3.2
while in the other cases, we will compare 
initial data $G_{ini}$ coming from $(t_D,x_D)$ 
with upstream data $G_{up}$ coming from $(t_U, x_U)$ (Region II), 
initial data $G_{ini}$ from $(t_U,x_U)$ 
with downstream data $G_{down}$ from $(t_D,x_D)$ (Region III) 
and finally upstream data $G_{up}$ from $(t_U,x_U)$ 
and downstream data $G_{down}$ from $(t_D,x_D)$ (Region IV). 
I feel that everything can be done exactly in the same way
and finally the idea is to use the inf-morphism from Aubin, Bayen, Saint-Pierre (2006)
$$ N(t,x) = \min \left\{ N_{ini}(t,x) \ , \ N_{up}(t,x) \ , \ N_{down}(t,x) \right\}.$$
}

\section{Numerical solution methods}
\label{sec:numerical}

In the following two subsections we formulate two numerical solution methods to find the global solution of the endogenous inflow problem in conservation law form (\ref{k-model}), (\ref{ephi}) and in HJ form (\ref{N-model}), (\ref{ephi}), respectively. The basic idea is that in discrete time, if the (endogenous) inflows are computed using the traffic states from the previous time step
{using an explicit-in-time numerical scheme}, 
then they become exogenous for the current time step and therefore the VT solution may be applied.

\subsection{Godunov's method}
This method has been traditionally used to solve the traffic
problem in conservation law form without inflows, and constitutes the basis of the well known Cell Transmission (CT) model \citep{Dag94}
{assuming the triangular Hamiltonian~\eqref{eq:triangFD}}. 
In this method, time and space are discretized in increments $\D t$ and $\D x=u \D t$, respectively, and we let:
\begin{equation}\label{}
k_i^j=k(j \D t,i \D x)
\end{equation}
be the numerical approximation of the density. The update scheme is the following discrete approximation of the conservation law (\ref{k-model}): 
\begin{equation}\label{ki}
\frac{k_{i+1}^{j}-k_{i}^{j}}{\D t}+\frac{q_{i}^{j+1}-q_{i}^{j}}{\D x}
=\phi(k_{i}^{j})
\end{equation}
The key to Godunov's method is the computation of the flow into cell  $i, q_{i}^{j}$, which are obtained by solving Riemann problems. Traditionally, inflows have been considered explicitly only in the update scheme (\ref{ki}) but not in the solution of the Riemann problems \citep{Lav10}. This implies that the computation of  $q_{i}^{j}$ corresponds to the original CT rule:
\begin{equation}\label{}
q_{i}^{j}=\min\{Q,u k_{i}^{j},(\kappa-k_{i}^{j+1})w\}, \hskip 3cm \mbox{(CT rule)}
\end{equation}

Here, we will compare the CT rule with the ERP rule (\ref{avgflow}), i.e. the flow based on extended Riemann problems. Both methods are first-order accurate since both are Godunov-type methods, and therefore the rate of convergence of both methods should be similar and roughly proportional to the mesh size. The main difference is in the magnitude of the error, where the ERP rule should be more accurate because the impacts of inflows are explicitly considered in the solution of Riemann problems. We illustrate this with the following example.
\\

\begin{figure}[tb]
\begin{center}
\includegraphics[width=.9\textwidth]{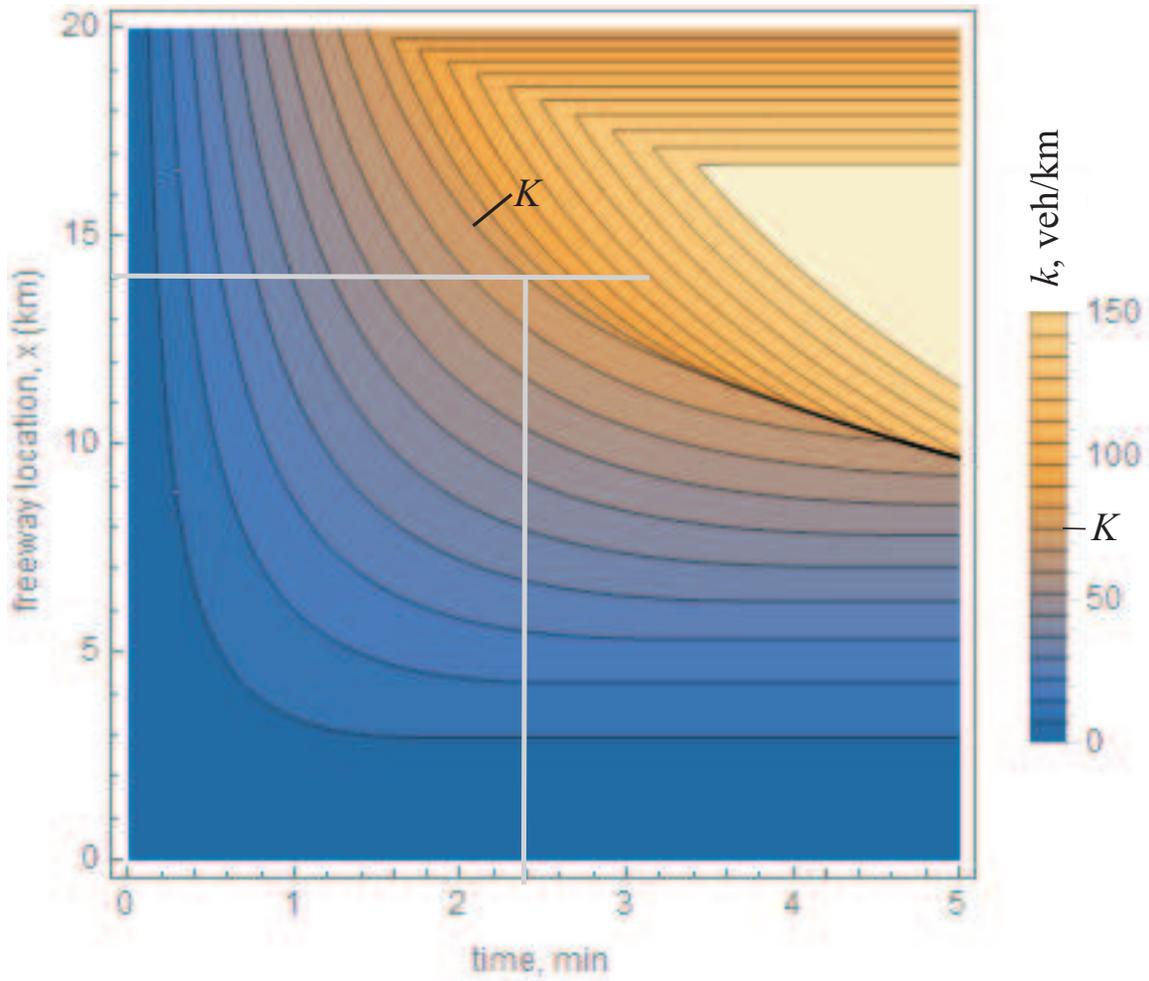}
\end{center}
\caption[]{Numerical ERP solution of  (\ref{exphi}) with $\D t=1$ s and  parameters: freeway length $L$ = 20 km, $w$ = 100 km/hr, $\kappa$ = 150 veh/km, $\theta$ = 1, $a = 0.5 Q/L$, $b = 0.3$.}
\label{f3}
\end{figure}

\emph{Example.} Consider an empty freeway at $t=0$ subject to an inflow linear in both $x$ and $k$; i.e.:
\begin{subequations}\label{exphi}
\begin{align}
	g(x)&=0, \\
	\phi(k)&=ax-buk, \hskip 1cm a,b>0.\label{inf}
\end{align}
\end{subequations}
Notice that \citep{Lav10} showed that linear inflow functions arise in the continuum approximation of the Newell-Daganzo merge model  \citep{New82,Dag94}, which accounts for the interactions between freeway and on-ramp demands. The particular coefficients of the linear function depend upon the state of the freeway and  on-ramps. In particular,  (\ref{inf}) corresponds to the case where both are in free-flow, $a x$  represent the inflow demand rate at $x$ and $b$ the exit probability per unit distance. Using  the method of characteristics \citep{Lav10} showed that the solution of  (\ref{k-model}), (\ref{exphi}) is :
\begin{equation}\label{eN}
k(t,x)=\frac{a}{b^2 u} \left(b x -1 +(1-b(x- t u))e^{-b t u}\right)
\end{equation}  
provided $k(t,x)\le K$. To get an idea of the solution for all densities, Fig.~\ref{f3} shows the  numerical ERP solution with $\D t=1$ s, with parameters given in its caption. Notice that similarly to \cite{Lav10} we chose $\theta=1$ to avoid the numerical errors intrinsic to Godunov's method and to focus on those caused by the treatment of the inflow.

\begin{figure}[tb]
\begin{center}
\includegraphics[width=\textwidth]{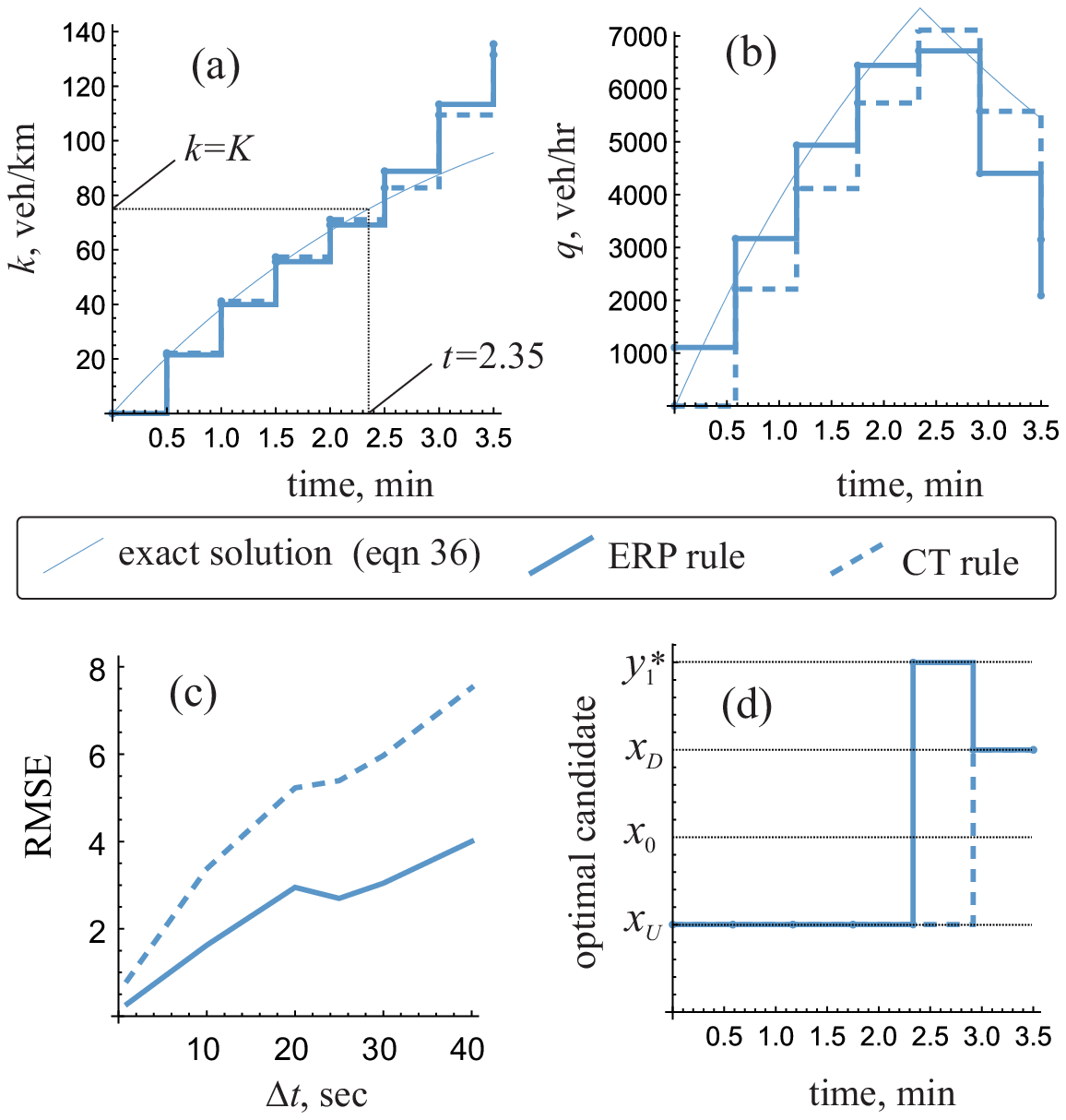}
\end{center}
\caption[]{Comparison of the CT rule with the ERP rule solution for the example in \S 4. (a) and (b): density and flow at $x=14$ km, (c)  candidate number that minimizes $f(y)$, (d) RMSE of each method for varying $\D t$.}
\label{f4}
\end{figure}

Fig.~\ref{f4}  compares the CT-rule with the ERP-rule numerical solution (with $\D t=40$ s) for this example \emph{vis-a-vis} the exact solution  (\ref{eN}). Parts (a) and (b) show the time evolution of the density and flow, respectively, at $x=14$ km obtained with each method. It can be seen that the main difference, as expected, is in the flow estimates, particularly at $t=0$ where the CT rule predicts zero flow since the freeway is empty and it does not consider inflows in its calculations. 

To assess the accuracy of each method, part (c) of  Fig.~\ref{f4} shows the density root-mean-squared error (RMSE) of each method with respect to (\ref{eN}) for varying $\D t$ and only until $t=2.35$ s, when the density exceeds the  critical density $K$; see Fig.~\ref{f4}a. It becomes apparent that  both methods converge to the right solution as $\D t\rightarrow 0$, but the accuracy of the proposed method outperforms the existing method by a factor of two for all values of $\D t$. 

Fig.~\ref{f4}d shows the optimal candidate that minimizes $f(y)$ at each time step of the numerical method, which is an element of the set $\BPP$, for the ERP rule, and of $\{x_U,x_0,x_D\}$ for the CT rule. It can be seen that both methods coincide except for the time step where the density approaches the critical density, in which case the proposed method finds the more accurate optimal candidate $y_1^*$.

\subsection{Variational networks}
\begin{figure}[tb]
\begin{center}
\includegraphics[width=.8\textwidth]{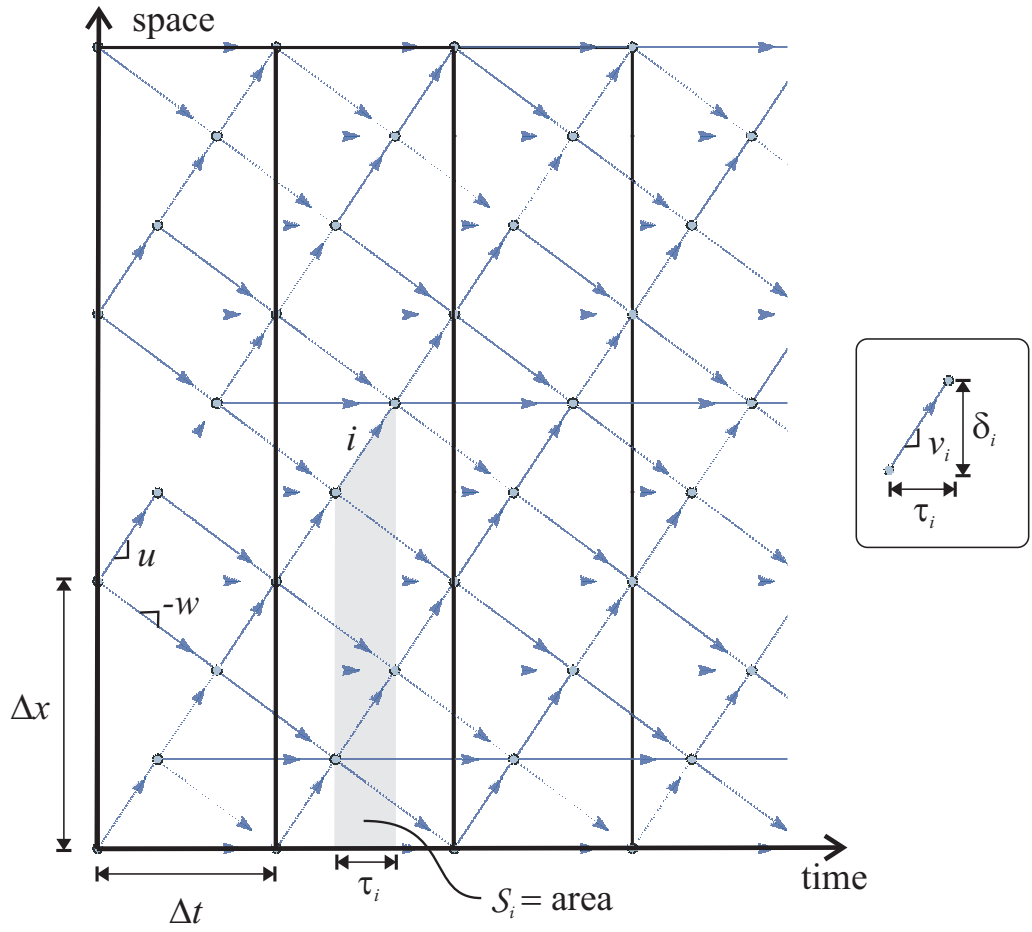}
\end{center}
\caption[]{Variational networks: time-space networks to solve the traffic problem with inflows using shortest paths.}
\label{f5}
\end{figure}

\cite{Dag05b} introduced time-space networks to solve the traffic problem without inflows in variational form using shortest paths.
{Notice that this is an application of Bellman's dynamic programming principle.}
Each link $i$ in these ``variational networks" is defined by its: (i) slope $v_i$: wave speed, (ii) cost $c_i$: maximum number of vehicles that can pass, (iii)  time length $\tau_i $, and (iv) distance length $\delta_i=\tau_i v_i $; see inset in Fig.~\ref{f5}.

Since the fundamental diagram is assumed triangular and the freeway homogeneous, there are only three wave speeds to be considered, $u, -w$ and $0$, and the corresponding passing rates  are given by
\begin{subequations}\label{NconstSol}
\begin{empheq}[left={\L(v_i)=\empheqlbrace}]{align}
     &w\kappa , &v_i&=-w\\
     &Q , &v_i&=0\\
     &0, &v_i&=u
\end{empheq}
\end{subequations}
Let  $J_i$ be the contribution of the $J$-integral in the cost of each link $i$. It corresponds to the (negative of) integral of $\phi(t,x)$ over the shaded region in Fig.~\ref{f5}, $\mathcal{S}_i$, and can be approximated by:
\begin{subequations}\label{}
\begin{align}
J_i&=- \tau_i\sum_{j\in \mathcal{S}_i}\delta_j a_j    ,\label{b}
\end{align}
\end{subequations}
where $a_j$ is the inflow associated with link $j$ and $j\in \mathcal{S}_i$ means all links that ``touch" area $\mathcal{S}_i$. Finally, the cost to be used in each link becomes:
\begin{equation}\label{}
c_i=\L(v_i)\tau_i + J_i.
\end{equation}

The advantage of this method is that it is free of numerical errors (when inflows are exogenous) but it may be cumbersome to implement unless $\theta$ is an integer. In that case, as illustrated in Fig.~\ref{f5} for $\theta=2$, the location of nodes align on a grid pattern. This allows defining a conventional grid with cell size $\D t, \D x$ where inflows may be assumed constant.  The other disadvantage is that merge models are typically expressed in terms of flows or densities rather than $N$ values, and therefore an additional computational layer has to be added.

\section{Other coordinates}
\label{sec:othercoord}

As pointed out in \cite{Lav13hj} there are two additional coordinate systems that provide alternative solution methods of the traffic flow problems without inflows. In space-Lagrangian coordinates  the quantity of interest is $X(t,n)$, the position of vehicle $n$ at time $t$; in time-Lagrangian coordinates one is interested in $T(n,x)$, the time vehicle~$n$ crosses location~$x$.

These representations correspond to the same surface in the three-dimensional space of vehicle number, time and distance, 
but expressed with respect to a different coordinate system. 
We briefly analyze these two alternatives 
{when considering Eulerian source terms}
and conclude that 
{even if a HJ equation is still valid in both cases,
one cannot expect to get a VT representation formula} 
even when inflows $\phi(t,x)$ are exogenous.

\subsection{Space-Lagrangian coordinates: X-models}
Let $s(t,n)$ be the spacing of vehicle $n$ at time $t$. To derive the X-model we multiply the N-model (\ref{N-modela}) by $s $ to get $sN_t-sH(k)=s\Phi(t,x)$. Noting that $s=-X_n$ and $sN_t=-X_nN_t=X_t$, this can be rewritten as: 
\begin{equation}\label{X-model}
	X_t-V(-X_n) =- X_n \Phi(t,X),
\end{equation}
where $V(s)=sH(1/s) $ is this spacing-speed fundamental diagram.  We conclude that (\ref{X-model})  is still a HJ PDE but does not admit a VT solution 
due to the term involving $X$. The corresponding conservation law can be obtained by taking the partial derivative with respect to $n$ of  (\ref{X-model}):
\begin{equation}\label{s-model}
	s_t+V(s)_n =-\phi(t,X)s^2-\Phi(t,X)s_n,
\end{equation}
Notice that \cite{femke} identified  (\ref{s-model}) but without the term $\Phi(t,X)s_n $ using a different approach, and used it to formulate a numerical solution method in the case of discrete inflows. 

\subsection{Time-Lagrangian coordinates: T-models}
Let $r=T_x$ and $h=1/q$ be the pace and the headway of vehicle $n$ at location $x$, and let $F(r)$ be the fundamental diagram in this case, i.e. $h=F(r)$. Here,  the T-model is simply $F(r)= 1/q$, where $q$ is given by  $q = N_t-\Phi(t,x)$, per  (\ref{flow}).  Noting that $N_t=1/T_n$ this can be rewritten as: 
\begin{equation}\label{T-model}
	T_n-\frac{F(T_x)}{1+\Phi(T,x)F(T_x)} =0,
\end{equation}
which, again, is still a HJ PDE but does not admit a VT solution
due to the term involving $T$. The corresponding conservation law can be obtained by taking the partial derivative with respect to $x$ of  (\ref{T-model}):
\begin{equation}\label{r-model}
	r_n-\frac{ F(r)_x+F(r)^2 \phi(T,x)}{(1+\Phi (T,x) F(r))^2}=0.
\end{equation}

To summarize, it becomes apparent that in Lagrangian and vehicle number-space coordinates the solution to our problem does not accept VT solutions and therefore becomes more difficult to solve. 

\section{Discussion}
\label{sec:disc}

We have shown in this paper that VT solutions to the traffic flow problem exist only in Eulerian coordinates when inflows are exogenous. In all other cases the Hamiltonian is a non-local function of the independent variable, and the corresponding variational may not be possible to formulate. Even in the simplest endogenous linear case (\ref{ephi}) the reader can appreciate the mathematical difficulties:  using $\Phi=\Phi(s,N(s,\xi(s)))$ in (\ref{Nsola})-(\ref{fxB00})  turns the  problem implicit in $N$, and therefore it is no longer a VT problem.

Improved numerical solution methods for the endogenous case were derived by taking advantage of this insight. In other fields, it appears that solving the extended Riemann problems explicitly considering the inflows has not been possible, and the only alternative has been to use high-resolution Riemann solvers \citep{Sch96,Lev98}. We have shown  that this is not the case in traffic flow, and that the ERP method presented here is indeed more accurate.

A streamlined version of the ERP method could be envisioned that drastically improves computation times with minimal impact in the quality of the solution. To see this, recall that Fig.~\ref{f4}d showed that the  candidate $y_1^*$ was optimal only during the time step where the density approaches the critical density. Based on the discussion following eqn.  (\ref{NconstSol}) it is reasonable to conjecture that candidates $y_1^*,y_2^*\ldots y_6^*$ in  (\ref{Nsolriem}) would be optimal only when a transition takes place. Therefore, the streamlined method would consider only the reduced set $\BPP=\{x_U,x_0,x_D\}$ in  (\ref{Nsolriem}), which would induce an error only in the time step where the transition occurs. This error can be made arbitrarily small by decreasing $\Delta t$. Notice that this does not mean that the continuum solution (\ref{Nsolriem}) can be streamlined in this way; this is only possible in discrete time where  $N$-values are updated at each time step.

The implications of our findings in the context of MFD analytical approximation methods considering turns as a continuum inflow are not encouraging. This is because inflows in this case would have to be  endogenous for the method to be meaningful, and in such case we have seen that there is no VT solution. This implies that the method of cuts--the only method used so far to provide analytical MFD  approximations--is no longer applicable. At the same time, however, a stochastic extension of the method of cuts proved successful in approximating a real-life MFD \citep{Lav15}. This would indicate that, at least in the context of the MFD, VT solutions still provide good approximations. Research in this topic is ongoing.

\section*{Acknowledgements}

This research was supported by NSF Grants 1055694 and 1301057. The author is grateful for the comments of two anonymous reviewers, which greatly improved the quality of this paper.

\section*{References}
\bibliographystyle{elsarticle-harv}
\bibliography{bibnew}

%\appendix
%\section{}
%This appendix 

\end{document}